\documentclass[11pt]{amsart}

\usepackage{latexsym,amssymb,amsmath,graphics,showkeys}
\usepackage[dvips]{graphicx}

\def\cc{{\mathcal C}}

\def\ss{{\mathcal S}}

\def\ie{i.e.\ }
 
\def\ffi{\varphi}
\def\eps{\varepsilon}
\def\dst{\displaystyle}

\def\C{{\mathbb{C}}}

\def\R{{\mathbb{R}}}

\def\Z{{\mathbb{Z}}}

\newcommand{\norm}[1]{{\left\|{#1}\right\|}}

\newcommand{\abs}[1]{{\left|{#1}\right|}}
\newcommand{\scal}[1]{{\left\langle{#1}\right\rangle}}

\newenvironment{notation}[1][]{\vskip1pt\noindent\rm\textit{Notation}\,:\ }{\rm\vskip1pt}
\newenvironment{remark}[1][]{\vskip1pt\noindent\rm\textit{Remark #1}\,:\ }{\rm\vskip1pt}
\newenvironment{definition}[1][]{\vskip3pt\noindent\sl\textbf{Definition.}\ }{\rm\vskip3pt}
\newenvironment{example}[1][]{\vskip1pt\noindent\rm\textit{Example\ #1}\,:\ }{\rm\vskip1pt}

\begin{document}

\title{Uncertainty principles for orthonormal bases}
\author{Philippe Jaming}
\address{Universit\'e d'Orl\'eans\\
Facult\'e des Sciences\\ 
MAPMO - F\'ed\'eration Denis Poisson\\ BP 6759\\ F 45067 Orl\'eans Cedex 2\\
France}
\email{Philippe.Jaming@univ-orleans.fr}

\begin{abstract}
In this survey, we present various forms of the uncertainty principle
(Hardy, Heisenberg, Benedicks). We further give a new interpretation
of the uncertainty principles as a statement about the time-frequency
localization of elements of an orthonormal basis, which improves
previous unpublished results of H. Shapiro.

Finally, we show that Benedicks' result implies that solutions of the Shr\"odinger
equation have some (appearently unnoticed) energy dissipation property.
\end{abstract}




\thanks{Research partially financed by: {\it European Commission}
Harmonic Analysis and Related Problems 2002-2006 IHP Network
(Contract Number: HPRN-CT-2001-00273 - HARP)
}
\maketitle

\section{Introduction}

The uncertainty principle is a ``metamathematical'' statement that asserts that
\begin{center}
\begin{minipage}{12cm}
{\sl a function and its Fourier transfrom can not both be sharply localized.}
\end{minipage}
\end{center}

There are various precise mathematical formulations of this general fact, the most well known being those of Hardy
and of Heisenberg-Pauli-Weil. It is one of the aims of this survey to present various statements that can be interpreted as
uncertainty principles. The reader may find many extensions of the results presented here, other manifestations of
the uncertainty principle as well as more references to the vast litterature in the surveys \cite{FS,BD}
and the book \cite{HJ}. Part of this survey has serious overlaps with these texts.

We start with results around Heisenberg's Uncertainty Principle and show how this can be proved using the spectral theory
of the Hermite Operator. We also give a not so well known extension that shows that Hermite Functions are ``successive
optimals'' of Heisenberg's Uncertainty Principle. We complete this section with a new result joint with A. Powell \cite{JP}
that shows that the elements of an orthonormal basis and their Fourier Transforms have their means and dispersions
that grow at least like those of the Hermite Basis. This gives a quantitative version of an unpublished result
by H. Shapiro (see \cite{Sh1,Sh2}).

In a second section, we consider the problem of the joint smallness of the support and spectrum 
(support of the Fourier Transform)
of a function. For instance, it is easy to show that a function and its Fourier Transform can not both have compact support.
It was shown by Benedicks \cite{Be} that the same is true when one asks the support and spectrum to be of finite measure.
Further, it was independently proved by Amrein and Berthier \cite{AB} that the $L^2$-norm of a function is controled
by the $L^2$-norm of the function outside a set of finite measure and the $L^2$-norm of its Fourier transform outside
an other set of finite measure. We then show that this property implies that the solutions of the Shr\"odinger equation
with initial data supported in a set of finite measure have a part of their energy dissipated outside any set of finite measure
(this is joint unpublished work with L. Grafakos). We complete the section by mentioning the closely related problem
of energy concentration in a compact set studied by Landau, Pollak and Slepian and conclude it with the local
uncertainty principle of Faris, Price and Sitaram.

We pursue the paper by measuring concentration by the speed of decay away from the mean.
We first sketch a new proof of Hardy's Theorem due to B. Demange \cite{De}. This roughly states that
a function and its Fourier Transform can not both decrease faster than gaussians and, if they decrease like
a gaussian, then they are actually gaussians. We conclude the section by asking whether all elements of an orthonormal
sequence as well as their Fourier Transforms be bounded by a fixed $L^2$-function. We show that this implies that the
sequence is finite and give bounds on the number of elements, thus improving a previous result of H. Shapiro
\cite{Sh1,Sh2}.

We conclude the paper by reformulating some uncertainty principles in terms of solutions
of the heat or the Shr\"odinger equation.

\section{Heisenberg's uncertainty principle}

In this section, we are going to prove Heisenberg's Uncertainty Principle and show that the Hermite Basis is the
best concentrated in the time-frequency plane when concentration is measured by dispersion. This is done using
the spectral theory of the Hermite Operator.

\subsection{Notations}\ \\
For sake of simplicity, all results will be stated on $\R$. Generalizations to $\R^n$ are either strateforward using tensorization
or available in the litterature.

For $f\in L^2(\R)$, let

--- $\mu(f)=\frac{1}{\norm{f}_2^2}\int_\R t|f(t)|^2\mbox{d}t$ its \emph{mean}.

--- $\Delta^2(f)=\int |t-\mu(f)|^2|f(t)|^2\mbox{d}t$ its \emph{variance}
and $\Delta(f)=\sqrt{\Delta^2(f)}$ its \emph{dispersion}.

\smallskip

For $f\in L^2(\R)\cap L^1(\R)$, the Fourier Transform is defined by
$$
\widehat{f}(\xi)=\int_\R f(x)e^{-2i\pi x\xi}\mbox{d}x
$$
and is then extended to all of $L^2(\R)$ in the usual way. The Inverse Fourier Transform of $f$ is denoted by $\check f$.

\smallskip

For $f\in\cc^2(\R)$ we define the Hermite operator as
$$
Hf(t)=-\frac{1}{4\pi^2}\frac{\mbox{d}}{\mbox{d}x}+t^2f(t).
$$

The Hermite Functions are then defined as
$$
h_k(t)=\frac{2^{1/4}}{\sqrt{k!}}\left(-\frac{1}{\sqrt{2\pi}}\right)^ke^{\pi t^2} \left(\frac{\mbox{d}}{\mbox{d}t}\right)^ke^{-2\pi t^2}.
$$
They have the following well known properties~:
\begin{enumerate}
\renewcommand{\theenumi}{\roman{enumi}}
\item they form an orthonormal basis in $L^2(\R)$,

\item they are eigenfunctions of the Fourier transform $\widehat{h_k}=i^{-k}h_k$,

\item they are eigenfunctions of the Hermite operator $\dst Hh_k=\frac{2k+1}{2\pi}h_k$.
\end{enumerate}

This allows to extend the Hermite Operator as a semi-bounded self-adjoint operator on $L^2(\R)$ via
$$
Hf=\dst\sum_{k=0}^{+\infty}\frac{2k+1}{2\pi}\scal{f,h_k}h_k.
$$
Finally, an easy computation shows that
$$
\scal{Hf,f}=\mu(f)^2\norm{f}^2+\Delta^2(f)+\mu(\widehat{f})^2\|\widehat{f}\|^2+\Delta^2(\widehat{f}).
$$
\subsection{Heisenberg's uncertainty principle}\ \\
It is now easy to give a proof of the following:

\medskip

\noindent{\bf Heisenberg-Pauli-Weil's Uncertainty principle.} {\sl For every $f\in L^2(\R)$, $\Delta(f)\Delta(\widehat{f})\geq\frac{1}{4\pi}\norm{f}_2^2$.
Moreover equality occurs only if there exists $C\in\C$, $\omega,a,\alpha\in\R$, such that}
$$
f(t)=ce^{2i\pi\omega t}e^{-\pi\alpha|t-a|^2}\quad a.e.
$$

\medskip

\begin{proof}
First, replacing $f$ by $f(t)=e^{2i\pi\omega t}f(t-a)$, it is enough to assume that $\mu(f)=\mu(\widehat{f})=0$.
Further, it is enough to prove that
$$
\Delta(f)^2+\Delta(\widehat{f})^2\geq \frac{1}{2\pi}\norm{f}_2^2
$$
with equality only for $f=c\,h_0$.
Indeed, it is enough to apply this to $f_\lambda(t)=\frac{1}{\sqrt{\lambda}}f(t/\lambda)$ and then minimize the left hand side
over $\lambda$.

But then, 
\begin{eqnarray*}
\Delta^2(f)+\Delta^2(\widehat{f})
&=&\scal{Hf,f}\\
&=&\sum_{k=0}^{+\infty}\frac{2k+1}{2\pi}|\scal{f,h_k}|^2\\
&\geq& \sum_{k=0}^{+\infty}\frac{2\times 0+1}{2\pi}|\scal{f,h_k}|^2=\frac{1}{2\pi}\norm{f}^2.
\end{eqnarray*}
Moreover, equality can only occur if $\scal{f,h_k}=0$ for all $k\not=0$, that is, for $f=ch_0$.
\end{proof}

\begin{remark} The above proof actually shows slightly more. If $f$ is orthogonal to $h_0,\ldots,h_{n-1}$, then
$$
\Delta^2(f)+\Delta^2(\widehat{f})\geq\frac{2n+1}{2\pi}\norm{f}^2
$$
with equality if and only if $f=ch_n$. This shows that the elements of the Hermite Basis are succesive optimizers of Heisenberg's Uncertainty Principle.
\end{remark}

\subsection{A quantitative form of Shapiro's Mean-Dispersion Theorem}\ \\
In an unpublished note, H. Shapiro \cite{Sh1,Sh2} asked what conditions should be put on sequences $(a_k)$, $(b_k)$, $(c_k)$, $(d_k)$ for
the existence of an orthonormal basis $(e_k)$ of $L^2(\R)$ such that
$$
\mu(e_k)=a_k,\ \mu(\widehat{e_k})=b_k,\ \Delta(e_k)=c_k,\ \Delta(\widehat{e_k})=d_k.
$$
Further he proved, using a compactness argument, that an orthonormal sequence such that all four
$\mu(e_k)$, $\mu(\widehat{e_k})$, $\Delta(e_k)$, $\Delta(\widehat{e_k})$ are bounded is necessarily finite.
These results were subsequently improved by A. Powell \cite{Po}
who proved that there is no orthonormal basis for which $\mu(e_k)$, $\Delta(e_k)$, $\Delta(\widehat{e_k})$ are bounded
and modified a construction of Bourgain to get an orthonormal basis for which $\mu(e_k)$, $\mu(\widehat{e_k})$, $\Delta(e_k)$
are bounded.

The question was motivated by the developement of Gabor and wavelet analysis. More precisely:
\begin{enumerate}
\item A \emph{Gabor basis} is a basis $\{g_{k,l}\}_{k,l\in\Z}$ of $L^2(\R)$ of the form $g_{k,l}(x)=e^{2i\pi kx}g(x-l)$.
It follows that $\mu(g_{k,l})\to\infty$ when $l\to\infty$ and $\mu(\widehat{g_{k,l}})\to\infty$ when $k\to\infty$.
Moreover, a theorem of Balian \cite{Ba} and Low \cite{Lo} asserts that if $(g_{k,l})_{k,l}$ is an orthonormal basis
then $\Delta(g)\Delta(\widehat{g})=+\infty$.

\item A \emph{Wavelet} is a basis $\{g_{k,l}\}_{k,l\in\Z}$  of $L^2(\R)$ of the form $g_{k,l}(x)=2^{k/2}g(2^kx-l)$.
Again $\mu(g_{k,l})\to\infty$ when $l\to\infty$, $\Delta(g_{k,l})\to\infty$ when $k\to-\infty$ and $\Delta(\widehat{g_{k,l}})\to\infty$ when $k\to+\infty$.
\end{enumerate}

We will now show the following quantitative form of Shapiro's result:

\smallskip

\noindent{\bf Theorem (Jaming-Powell \cite{JP}).} {\sl
Let $\{e_k\}_{k\geq0}$ be an orthonormal sequence in $L^2(\R)$.
Then for all $n\geq 0$,
$$
\sum_{k=0}^n \left( \Delta^2(e_k) + \Delta^2(\widehat{e_k}) + 
|\mu(e_k)|^2 + |\mu(\widehat{e_k})|^2 \right) \geq \frac{(n+1)^2}{2\pi}.
$$
If equality holds for all $n\leq n_0$, then
for $k=0,\ldots,n_0$, $e_k=c_kh_k$, $|c_k|=1$.}

\smallskip

As a corollary, we immediatly get the following:

\smallskip

\noindent{\bf Corollary (Jaming-Powell \cite{JP}).} {\sl
Let $\{e_k\}_{k\geq0}$ be an orthonormal sequence in $L^2(\R)$
such that $\Delta(e_k)$, $\Delta(\widehat{e_k})$, $|\mu(e_k)|$, $|\mu(\widehat{e_k})|$ are all
$\leq C$
then the sequence has at most $8\pi C^2$ elements.}

\smallskip

The proof is a direct application of the Rayleigh-Ritz technique which we now recall:

\smallskip

\noindent{\bf Theorem (Rayleigh-Ritz).} {\sl
Let $H$ be a semi-bounded self-adjoint operator with domain $D(H)$. Define
$$
\lambda_k (H) = \sup_{\varphi_0, \cdots, \varphi_{k-1}} \ \ 
\inf_{\psi \in [\varphi_0, \cdots, \varphi_{k-1}]^{\perp},||\psi||=1, 
\psi \in D(H)} 
\langle H \psi, \psi \rangle.
$$
Let
$V$ be a $n+1$ dimensional
subspace $V\subset D(H)$ and let $P_V$ be the orthogonal projection on $V$.
Define $H_V = P_VHP_V$, and consider $H_V$ as an operator $\widetilde{H_V}:\ V\to V$. Let $\mu_0 \leq \mu_1 \leq \cdots \leq \mu_n$ be the eigenvalues of 
$\widetilde{H_V}$. Then}
$$
\lambda_k (H) \leq \mu_k, \ \ \ k=0, \cdots, n.
$$

\smallskip

Taking the trace of $\widetilde{H_V}$, we then get that, for $(\varphi_k)$ orthonormal and in $D(H)$, 
$$
\sum_{k=0}^n \lambda_k(H)  \leq \sum_{k=0}^n
\langle H\varphi_k, \varphi_k \rangle.
$$
If one of the $\varphi_k$'s is not in $D(H)$, then the right-hand side is infinite and this is trivial.
To conclude the proof of the theorem, it is thus enough to take
$H$ the Hermite operator and $\varphi_k=e_k$.\hfill$\Box$

\section{Qualitative uncertainty principles}

\subsection{Annihilating pairs}
\ 
\begin{notation}
We will use the following notations. All sets considered will be measurable and the measure of a set $S$ will be denoted by $|S|$.
We denote by $\chi_S$ the characteristic function of $S$, that is $\chi_S(x)=1$ if $x\in S$ and $\chi_S(x)=0$ otherwise.

For $S,\Sigma\subset\R$, let $P_S,Q_\Sigma$ be the projections on $L^2(\R)$ defined by $P_Sf=\chi_S f$ and
$Q_\Sigma f=(\chi_\Sigma\widehat{f})\check{\ }$. Let $P^\bot_S=I-P_S=P_{\R\setminus S}$ and $Q^\bot_\Sigma=I-Q_\Sigma$.

Finally, let $L^2(S)$ be the space of all functions in $L^2(\R)$ with support in $S$ and $L^2_\Sigma$ be the set of all functions
whose Fourier Transform has support in $\Sigma$.
\end{notation}

A way of measuring the concentration of a function in the time-frequency plane is by asking that 
the function and its Fourier Transform have both small support.
For instance, it is well known that a function and its Fourier Transform can not both
have compact supported. Indeed if a non zero function has compact support, then its Fourier tTransform
is analytic and thus can not have compact support. This leads naturally to the following property:

\begin{definition}
Let $S,\Sigma$ be two measurable subsets of $\R^d$. Then 

\begin{itemize}
\item[---] $(S,\Sigma)$ is a (weak) \emph{annihilating pair} if, 
$$
\mbox{supp}\,f\subset S\quad\mbox{and}\quad\mbox{supp}\,\widehat{f}\subset\Sigma
$$
implies $f=0$.

\item[---] 
$(S,\Sigma)$ is called a \emph{strong annihilating pair} if there exists $C=C(S,\Sigma)$
such that
\begin{equation}
\label{eq:sa1}
\norm{f}_{L^2(\R^d)}\leq C\bigl(\norm{f}_{L^2(\R^d\setminus S)}+\|\widehat{f}\|_{L^2(\R^d\setminus\Sigma)}\bigr).
\end{equation}
\end{itemize}
\end{definition}

To prove that a pair $(S,\Sigma)$ is a strong annihilating pair, one usually shows that there exists a constant 
$D=D(S,\Sigma)$ such that, for all functions $f\in L^2(\R^d)$ whose Fourier Transform is
supported in $\Sigma$,
\begin{equation}
\label{eq:sa2}
\norm{f}_{L^2(\R^d\setminus\Sigma)}\geq D\norm{f}_{L^2(\R^d)}.
\end{equation}
Moreover, the best constants are related by $\frac{1}{D(S,\Sigma)}\leq C(S,\Sigma)\leq\frac{1}{D(S,\Sigma)}+1$.

As seen as above, if $S,\Sigma$ are compact, then they form an annihilating pair. More generally, if 
$S,\Sigma$ are sets of finite measure, then they form an annihilating pair. 
This is a particular case of the following:

\smallskip

\noindent{\bf Theorem (Benedicks \cite{Be}).} {\sl Assume that, for almost every $x\in(0,1)$, the set $(x+\Z)\cap S$
is finite and, for almost every $\xi\in(0,1)$,
the set $(\xi+\Z)\cap\Sigma$ is finite. Then the pair $(S,\Sigma)$ is a weak annihilating pair.}

\begin{proof} The main tool is the Poisson Summation Formula, which implies that
\begin{equation}
\label{eq:poisson}
e^{-2i\pi x\xi}\sum_{j\in\Z}f(x+j)e^{2i\pi j\xi}=\sum_{k\in\Z}\widehat{f}(\xi+k)e^{2i\pi kx}.
\end{equation}

Assume first that there exists a set $E\subset(0,1)$ of positive measure such that, for all $x\in E$, the set $(x+\Z)\cap S$ is empty.
Then the left hand side of (\ref{eq:poisson}) vanishes on $E$. Since the right hand side is a trigonometric polynomial in the $x$ variable, for almost every $\xi$,
 it is identically $0$. In particular
 $$
\norm{f}_2^2= \int_0^1 \sum_{k\in\Z}|\widehat{f}(\xi+k)|^2=0
$$
and the conclusion follows in this case.

Let us now prove the general case. It is easy to find a set $E\subset(0,1)$ of positive measure and an integer $N>0$
such that $(x+N\Z)\cap S$ is empty. We will use the previous case after a change of scale. We are thus
lead to proving that, for almost all $\xi\in(0,N^{-1})$, the set $(\xi+N^{-1}\Z)\cap\Sigma$ is finite,
which is done by writing it as a finite union of sets $(\xi_j+\Z)\cap\Sigma$ where $\xi_j=\xi+\frac{j}{N}$, $j=0,\ldots,N$.
\end{proof}

It has been proved by Amrein-Berthier \cite{AB} that pairs of sets of finite measure are actually strongly annihilating.
As noticed in \cite{BD}, this can be deduced from Benedick's theorem.

\smallskip

\noindent{\bf Theorem (Amrein-Berthier \cite{AB}).} {\sl Let $S,\Sigma$ be sets of finite measure, then
$(S,\Sigma)$ is weak annihilating pair.}

\begin{proof} Assume there is no such constant $C$. Then there exists a sequence $f_n\in L^2(\R)$ of norm $1$ and with spectrum 
$\Sigma$ such that $f_n\chi_{\R\setminus S}$ converges to $0$. Moreover, we may assume that 
$\widehat{f_n}$ is weakly convergent in $L^2(\R)$ with some limit $\widehat{f}$.
As $f_n(x)$ is the scalar product of $e^{2i\pi x\xi}\chi_\Sigma$ and $\widehat{f_n}$, it follows that $f_n$ converges to $f$.
Finally, as $|f_n|$ is bounded by $|\Sigma|^{1/2}$, we may apply Lebesgue's Theorem, thus $f_n\chi_S$ converges to $f$
in $L^2(\R)$ and the limit $f$ has norm $1$. But the function $f$ has support in $S$ and spectrum in $\Sigma$ so by 
Benedick's Theorem,
it is $0$, a contradiction.
\end{proof}

It is much more difficult to have an estimate of the constant in Amrein-Berthier's Theorem. Let us mention a theorem of
Nazarov \cite{Na} which states that one can take $C=c_0e^{c_1|S||\Sigma|}$. The same problem in higher dimension is not yet solved in a completely satisfactory way.

For particular $\Sigma$, the range of possible $S$ may be widened. Let us mention the following theorem, which improves
estimates of Logvinenko-Sereda \cite{LS} and Paneah \cite{Pa1,Pa2}. 

\begin{definition} A set $E$ is $\gamma$-thick at scale $a>1$ if, for all $x\in\R$,
$$
|E\cap[x-a,x+a]|\geq 2\gamma a.
$$
\end{definition}

\noindent{\bf Theorem (Kovrizhkin \cite{Ko}).} {\sl There exists a constant $C$ such that, for all functions $f\in L^2(\R)$
with spectrum in $[-1,1]$ and for every set $E$ which is $\gamma$-thick at scale $a>1$, one has the inequality}
$$
\norm{f}^2_2\le \left(\frac{C}{\gamma}\right)^{Ca}\int_{E}|f(x)|^2\mbox{d}x.
$$

\smallskip

\subsection{Prolate spheroidal wave functions}

A closely related question is that of concentration. This has been extensively studied by Landau, Pollak and Slepian
\cite{LP1,LP2,SP,S1}, see also \cite{S2}. More precisely, they studied the set of all functions $f\in L^2(\R)$
such that, for some fixed $T,\Omega,\eps>0$,
\begin{equation}
\label{eq:prolate}
\int_{|t|>T}|f(t)|^2\mbox{d}t\leq\varepsilon^2\norm{f}_2^2\quad,\quad
\int_{|\xi|>T}|\widehat{f}(\xi)|^2\mbox{d}\xi\leq\varepsilon^2\norm{f}_2^2.
\end{equation}
In particular, they proved that there is an orthonormal basis $\{\psi_n\}_{n\geq 0}$ of the space
$$
PW_\Omega=\{f\in L^2(\R)\,:\ \mbox{supp}\,\widehat{f}\subset[-\Omega,\Omega]\}
$$ 
formed of eigenfunctions  of
the differential operator:
$$
L=(T^2-x^2)\left(\frac{\mbox{d}}{\mbox{d}x}\right)^2-2x\frac{\mbox{d}}{\mbox{d}x}-\frac{\Omega^2}{T^2}x^2
$$
such that the following theorem holds:

\smallskip

\noindent{\bf Theorem (Landau-Pollak \cite{LP2})} {\sl Let $d=\lfloor 4T\Omega\rfloor+1$. Then,
for every $f\in L^2$, s.t. $(*)$ holds
$$
\norm{f-\mathbb{P}_df}_{L^2(\R)}^2\leq 49\varepsilon^2\norm{f}_2^2,
$$
where $\mathbb{P}_d$ is the orthogonal projection on the span of $\psi_0,\ldots,\psi_{d-1}$.}

\smallskip

The functions $\psi_k$ are called \emph{prolate spheroidal wave functions} and are linked to the wave
equation in ``prolate spheroidal'' coordinates.

This theorem is intimately linked to the Shannon Sampling Theorem. Its heuristics is the following.
Assume that $f\in PW_\Omega$ then, according to Shannon's sampling theorem, $f$ can be reconstructed from its samples
in the following way~:
$$
f(x)=\sum_{k\in\Z}f\left(\frac{k}{2\Omega}\right)\mbox{sinc}\,\pi\left(x-\frac{k}{2\Omega}\right).
$$
But now, if $f$ is essentially $0$ outside $[-T,T]$, then $f$ is essentially reconstructed from its $4T\Omega$ samples
at $\{\frac{k}{2\Omega},-2T\Omega\leq k\leq 2T\Omega\}$.
That is, the space of functions that are almost supported in $[-T,T]$ and have spectrum in $[-\Omega,\Omega]$ has dimension $d\simeq4T\Omega$.
The above theorem gives a precise formulation of this.

\subsection{Local uncertainty principles}

The following theorem is due to Faris\cite{Fa} when $\alpha=1$ and Price and Sitaram \cite{Pr,PS} in the general case.

\smallskip

\noindent{\bf Theorem.} {\sl If $0<\alpha<d/2$, there is a constant $K=K(\alpha,d)$ such that for all $f\in L^2(\R^d)$ and all measurable
set $E\subset\R^d$,}
\begin{equation}
\label{eq:faris}
\int_E|\widehat{f}(\xi)|^2\mbox{d}\xi\leq K|E|^{2\alpha/d}\norm{|x|^\alpha f}_2^2.
\end{equation}

\smallskip

\begin{proof}
Let $\chi_r$ be the characteristic function of the ball $\{x\,: |x|<r\}$ and $\tilde\chi_r=1-\chi_r$. Then, for $r>0$, write
$$
\left(\int_E|\widehat{f}(\xi)|^2\mbox{d}\xi\right)^{1/2}=\|\widehat{f}\chi_E\|_2\leq\|\widehat{f\chi_r}\chi_E\|_2
+\|\widehat{f\tilde\chi_r}\chi_E\|_2
\leq |E|^{1/2}\|\widehat{f\chi_r}\|_\infty+\|\widehat{f\tilde\chi_r}\|_2.
$$
Now
$$
\|\widehat{f\chi_r}\|_\infty\leq\|f\chi_r\|_1\leq\|\,|x|^{-\alpha}\chi_r\,\|_2\|\,|x|^\alpha f\,\|_2\leq C_\alpha r^{d/2-\alpha}\|\,|x|^\alpha f\,\|_2
$$
and
$$
\|\widehat{f\tilde\chi_r}\|_2\leq \|\,|x|^{-\alpha}\tilde\chi_r\,\|_\infty\|\,|x|^\alpha f\,\|_2\leq r^{-\alpha}\|\,|x|^\alpha f\,\|_2,
$$
so
$$
\left(\int_E|\widehat{f}(\xi)|^2\mbox{d}\xi\right)^{1/2}\leq \bigl(C_\alpha |E|^{1/2}r^{d/2-\alpha}+r^{-\alpha}\bigr)\|\,|x|^\alpha f\,\|_2.
$$
The desired result is obtained by minimizing the right hand side of that inequality over $r>0$.
\end{proof}

An easy computation shows that this proof gives 
$$
K=\frac{(d+2\alpha)^2}{(2\alpha)^{4\alpha/d}}(d-2\alpha)^{2\alpha/d-2}.
$$

A more interresting remark is the following. First, note that if we exchange $f$ and $\widehat{f}$ in (\ref{eq:faris}), 
then the right hand side becomes, $\norm{\Delta^{\alpha/2} f}_2$, using Parseval.
We thus get that
$$
\sup_{E\subset\R^d,\ 0<|E|<\infty}\frac{1}{|E|^{\alpha/d}}\left(\int_E|f(x)|^2\mbox{d}x\right)^{1/2}
\leq K\norm{\Delta^{\alpha/2} f}_2.
$$
The left hand side is known to be an equivalent norm of the Lorentz-space $L^{p_\alpha,\infty}$ where $\dst p_\alpha=\frac{2d}{d-2\alpha}$~:
$$
\bigl\|f\bigr\|_{L^{p_\alpha,\infty}}\leq K\bigl\|\Delta^{\alpha/2} f\bigr\|_2\quad\mbox{ for all }f\in L^2(\R).
$$

\section{Fast decrease conditions}

We will now measure concentration by asking that a function has fast decrease away from its mean.
A typical very concentrated function is the gaussian $e^{-\pi ax^2}$ and the parameter $a$ measures
concentration.

\subsection{Hardy's uncertainty principle}\ \\[3pt]
\noindent{\bf Theorem (Hardy \cite{Ha}).} {\sl Let $f\in L^2(\R)$. Assume that
\begin{enumerate}
\item for almost all $x\in\R$, $|f(x)|\leq C(1+|x|)^Ne^{-\pi a|x|^2}$,

\item for almost all $x\in\R$, $|\widehat{f}(\xi)|\leq C(1+|\xi|)^Ne^{-\pi b|\xi|^2}$.
\end{enumerate}
Then, if $ab>1$, $f=0$ and
if $ab=1$ $f(x)=P(x)e^{-\pi a|x|^2}$, $P$ polynomial of degree at most $N$.}

\smallskip

Note that this theorem may be interpreted as follows~: the space of functions satisfying conditions (1) and (2) of Hardy's Theorem is $N+1$-dimensional.
There are various generalizations of this theorem. The following is an improvement of a result of
Beurling whose proof was lost until a new prove was given by H\"ormander. H\"ormander's result was only in dimension $1$
and with $N=0$ so that it did not cover the equality case in Hardy's Theorem.

\smallskip

\noindent{\bf Theorem (Bonami-Jaming-Demange \cite{BDJ}).} {\sl Let $f\in L^2(\R^d)$. Assume that
$$
\iint_{\R^d\times\R^d}|f(x)||\widehat{f}(\xi)|e^{2\pi\abs{\scal{x,\xi}}}\frac{\mbox{d}x\mbox{d}\xi}{(1+|x|+|\xi|)^N}
$$
then $f=P(x)e^{-\scal{Ax,x}}$, where $A$ is positive definite and $P$ is a polynomial of degree $<\frac{N-d}{2}$.}

\smallskip

This theorem was further generalized to distributions by B. Demange. Its only disadvantage is that it does not give
information on the degree of the polynomial, but in practice, this can be easily obtained.

\smallskip

\noindent{\bf Theorem (Demange).} {\sl Let $f\in\mathcal{S}'(\R^d)$. Assume that
$e^{\pi|x|^2}f\in\mathcal{S}'$ and $e^{\pi|\xi|^2}\widehat{f}\in\mathcal{S}'$ then 
$f=P(x)e^{-\scal{Ax,x}}$, where $A$ is positive definite and $P$ is a polynomial.}

\begin{proof}[Sketch of proof for $d=1$]
The key tool here is to use the Bargmann Transform which transforms the distribution $f$ into an entire function of order 2 via
$$
F(z)=\mathcal{B}(f)(z):=e^{\frac{\pi}{2}z^2}\scal{f(x),e^{-\pi(x-z)^2}}.
$$
Note that $F(-iz)=\mathcal{B}(f)(-iz)=\mathcal{B}(\widehat{f})(z)$.
Using the fact that a Schwartz distribution has finite order, we get from $e^{\pi|x|^2}f\in\mathcal{S}'$ that there exists $C,N$ such that
$$
|F(z)|\leq C(1+|z|)^Ne^{\frac{\pi}{2}|\mathrm{Im}z|^2}
$$
while, from $e^{\pi|\xi|^2}\widehat{f}\in\mathcal{S}'$ and the previous formula, we get that
$$
|F(z)|\leq C(1+|z|)^Ne^{\frac{\pi}{2}|\mathrm{Re}z|^2}.
$$
But then, Phragm\`en-Lindel\"of's Principle implies that $|F(z)|\leq C(1+|z|)^N$ which, with Liouville's Theorem, implies that
$F$ is a polynomial. Inverting the Bargmann Transform then shows that $f$ is of the desired form.
\end{proof}

\subsection{The umbrella theorem}

We will now prove that the elements of an orthonormal basis and their Fourier Transforms can not all be bounded
by fixed functions:

\smallskip

\noindent{\bf Theorem (Jaming-Powell).} {\sl Let $\varphi,\psi$ be two \emph{fixed} functions in $L^2(\R)$.
Then there exists $N=N(\varphi,\psi)$ such that, if $(e_k)_{k\in I}\subset L^2(\R)$ is orthonormal and
$$
|e_k|\leq \varphi\quad,\quad |\widehat{e}_k|\leq \psi
$$
then $\#I\leq N$.}

\smallskip

This theorem generalizes an earlier result of H. Shapiro \cite{Sh1,Sh2}
that showed, using a compactness argument, that the sequence is finite, but with no quantitative estimate.

\begin{proof}
Let $M=\max(\norm{\varphi},\norm{\psi})$, and $0<\varepsilon<\frac{1}{50M}$. Define
$$
C_\varphi(\varepsilon)=\inf\left\{T\in\R\,:\ \int_{|t|>T}|\varphi(t)|^2\leq\varepsilon^2\norm{\varphi}_2^2
\right\}.
$$
Define $C_\psi(\varepsilon)$ accordingly.
For $T>\max\bigl(C_\varphi(\varepsilon),C_\psi(\varepsilon)\bigr)$, for all $n\in I$,
$$
\int_{|t|>T}|e_n(t)|^2\mbox{d}t\leq \varepsilon^2\norm{\varphi}_2^2\quad,\quad
\int_{|\xi|>T}|\widehat{e_n}(\xi)|^2\mbox{d}\xi\leq \varepsilon^2\norm{\psi}_2^2.
$$
Let $d=\lfloor 4T^2\rfloor+1$, $\eta=7M\varepsilon$ and let
$\psi_n$ be the associated prolate sphero\"idals,
$$
\mbox{for }n\in I,\quad\norm{e_n-\mathbb{P}_de_n}\leq\eta.
$$

Now write $\dst\mathbb{P}_de_n=\sum_{k=0}^{d-1}a_{n,k+1}\psi_k$ {\it i.e.} $a_n:=(a_{n,k})\in\C^d$. Then
$1-\eta\leq\norm{a_{n,k}}=\norm{\mathbb{P}_de_n}\leq 1$ and
\begin{eqnarray*}
\scal{a_{n,k},a_{m,k}}&=&\scal{\mathbb{P}_de_n,\mathbb{P}_de_m}\\
&=&\scal{\mathbb{P}_de_n-e_n+e_n,\mathbb{P}_de_m-e_m+e_m}\\
&=&\scal{\mathbb{P}_de_n-e_n,\mathbb{P}_de_m-e_m}+\scal{e_n,\mathbb{P}_de_m-e_m}+\scal{\mathbb{P}_de_n-e_n,e_m}\\
&=&\scal{\mathbb{P}_de_n-e_n,\mathbb{P}_de_m-e_m}+\scal{e_n-\mathbb{P}_de_n,\mathbb{P}_de_m-e_m}+
\scal{\mathbb{P}_de_n-e_n,e_m-\mathbb{P}_de_m}\\
&=&\scal{\mathbb{P}_de_n-e_n,e_m-\mathbb{P}_de_m}\leq \eta^2.
\end{eqnarray*}
So $b_n=\dst\frac{a_n}{\norm{a_n}}$ are vectors on the unit sphere such that the angles 
$|\scal{b_n,b_m}|\leq\frac{\eta^2}{1-\eta^2}$.
A set of such vectors is called a \emph{spherical code} in $\C^d$ and is known to be finite.
\end{proof}

There are various bounds on the number of spherical codes.
For instance, one may identify $\C^d$ with $\R^{2d}$ in the standard way. One then gets that
if $|\scal{b_n,b_m}_{\C^d}|\leq\alpha$ then $\scal{b_n,b_m}_{\R^{2d}}\in[-\alpha,\alpha]$
\ie $(b_n)$ is a \emph{$[-\alpha,\alpha]$-spherical code} in $\R^{2d}$.
If we call $N^{2d}(\alpha)$ the maximal number of elements of a $[-\alpha,\alpha]$-spherical codes in $\R^{2d}$,
then the following bounds can be obtained
\begin{itemize}
\item[---] if $\alpha<\frac{1}{2d}$ then $N^{2d}(\alpha)=2d$.\\
This is a trivial bound as the code is then linearly independent.

\item[---]  $N^{2d}(\alpha)\leq \left(\frac{2-\alpha}{1-\alpha}\right)^{2d}$ .\\
This can be obtained by a volume counting argument.

\item[---] if $\alpha<\frac{1}{\sqrt{2d}}$ then
$N^{2d}(\alpha)\leq2\frac{1-\alpha^2}{1-2\alpha^2d}d$,
has been obtained by Delsarte-Goethals-Siedel \cite{DGS}.
\end{itemize}

If $\ffi,\psi$ are explicitely given, one may then optimize the choice of $T$ and $\eps$ in the previous proof to get
more explicit estimates on $\#I$.

\begin{example}\ 
\begin{enumerate}
\item Take $\varphi=\psi=Ce^{-\pi a x^2}$ then
$$\dst
\#I\leq 2+\frac{2}{\pi a}\max\left(1+\ln\frac{C}{a^{5/2}},\pi+\frac{1}{2}\ln\frac{C^2}{2a}\right).
$$
In this case, only the trivial bound is needed.
\item Take $\varphi=\psi=\frac{C}{(1+|x|)^p}$ then, using the volume counting bound, one gets
$$
\#I\leq N\leq 3^{8\left(\frac{20C}{\sqrt{2p-1}}\right)^{\frac{4}{2p-1}}}.
$$
Moreover,
$$
\#I\leq \left\{\begin{matrix}
4\left(\frac{800C^2}{2p-1}\right)^{\frac{2}{2p-3}}&\mbox{ if }p>3/2\\
32\left(20C\right)^{\frac{2}{2p-3}}&\mbox{ if }1<p\leq 3/2\\
\end{matrix}\right.
$$
where we used the trivial bound for the first estimate and Delsarte {\it et al}'s bound for the second.
\end{enumerate}
\end{example}

\section{Reformulation of uncertainty principles as properties of solutions of some
classical PDE's}

\subsection{The Free Heat Equation}
Let us recall that the solution of the Free Heat Equation
\begin{equation}
\label{eq:heat}
\begin{cases}
-\partial_t v+\frac{1}{4\pi}\partial_x^2 v=0\\
v(x,0)=v_0(x)\\
\end{cases}
\end{equation}
with initial data $v_0\in L^2(\R)\cup L^1(\R)$ has solution
$$
v(x,t)=\int_{\R}e^{-\pi\xi^2t+2i\pi x\xi}\widehat{v_0}(\xi)\mbox{d}\xi
=\bigl(e^{-\pi\xi^2t}\widehat{v_0}\bigr)\check\ .
$$
We may then reformulate some uncertainty principles as follows.

\smallskip

\noindent{\bf Hardy's Uncertainty Principle}. {\sl Let $v_0\in L^1(\R)$ be non zero
and let $v$ be a solution of (\ref{eq:heat}) with initial data $v_0$. Then for every $t>0$
and for every $N>0$,
$\dst x\mapsto \frac{e^{\pi x^2/t}}{(1+|x|)^N}v(x,t)$ is unbounded.}

\smallskip

Indeed, as $e^{-\pi\xi^2t}|\widehat{v_0}(\xi)|\leq\norm{v_0}_{L^1(\R)}e^{-\pi\xi^2t}$, if
$|v(x,t)|\leq C(1+|x|)^Ne^{-\pi x^2/t}$ then, from Hardy's Theorem, $e^{-\pi\xi^2t}|\widehat{v_0}(\xi)|=Ce^{-\pi\xi^2t}$,
a contradiction.

Note that Demange's Theorem shows that, if $v_0\in\ss'(\R)$ and if $x\mapsto e^{\pi x^2/t}v(x,t)\in\ss'(\R)$
then $v_0$ is a Hermite Function, an so if then $v$.
\medskip

\noindent{\bf Annihilating pairs.} Let $(S,\Sigma)$ be an annihilating pair.
Let $v_0\in L^2(\R)$ and assume that $\widehat{v_0}$ is supported in $\Sigma$.
Let $v$ be a solution of (\ref{eq:shr}) with initial data $v_0$ and assume that,
for some $t>0$, $v(\cdot,t)$ is supported in $S$, then $v_0=0$.

Moreover, if $(S,\Sigma)$ is a strong annihilating pair with constant $D(S,\Sigma)$
in (\ref{eq:sa2}), then, for all $t>0$,
$$
\norm{v(x,t)}_{L^2(\R\setminus S)}\geq D\norm{v(x,t)}_{L^2(\R)}=D(S,\Sigma)\norm{e^{-\pi\xi^2 t}\widehat{v_0}(\xi)}_{L^2(\Sigma)}.
$$
For instance, if $\Sigma=[-a,a]$, we thus get the lower bound
$$
\norm{v(x,t)}_{L^2(\R\setminus S)}\geq D(S,\Sigma)e^{-\pi a^2t}\norm{v_0}_{L^2(\R)}.
$$

\subsection{The Free Shr\"odinger Equation}
Let us recall that the solution of the Free Shr\"odinger Equation
\begin{equation}
\label{eq:shr}
\begin{cases}
i\partial_t v+\frac{1}{4\pi}\partial_x^2 v=0\\
v(x,0)=v_0(x)\\
\end{cases}
\end{equation}
with initial data $v_0\in L^2(\R)$ has solution
$$
v(x,t)=\int_{\R}e^{-i\pi\xi^2t+2i\pi x\xi}\widehat{v_0}(\xi)\mbox{d}\xi
=\bigl(e^{-i\pi\xi^2t}\widehat{v_0}\bigr)\check\ .
$$
We may then reformulate some uncertainty principles as follows.

\medskip

\noindent{\bf Hardy's Uncertainty Principle}. {\sl Let $v_0\in L^2(\R)$ and assume that
$|\widehat{v_0}(\xi)|\leq Ce^{-\pi a\xi^2}$.
Let $v(x,t)$ is a solution of (\ref{eq:shr}) with initial data $v_0$ and assume that,
for some $t>0$, there exists $c=c(t)$ and $b=b(t)$ such that
$$
|v(x,t)|\leq ce^{-\pi b\xi^2}.
$$
If $b>1/a$, then $v_0=0$.
If $b=1/a$, then $v_0(x)=Ce^{-\pi x^2/(a-it_0)}$.
for some $C\in\C$.}

\smallskip

The first assertion is a direct consequence of Hardy's theorem. For the second,
Hardy's theorem implies that $v(x,t_0)=ce^{-\pi \xi^2/a}$, thus
$\widehat{v_0}(\xi)=ce^{i\pi\xi^2t_0}e^{-\pi a\xi^2}=ce^{-\pi(a-it_0)\xi^2}$.
Inverting the Fourier transform, we get the assertion.
Note that, in this case
$$
|v(x,t)|=\frac{C}{\sqrt{a^2+(t-t_0)^2}}\exp\left(-\frac{\pi a}{a^2+(t-t_0)^2}x^2\right).
$$

\medskip

\noindent{\bf Annihilating pairs.} Let $(S,\Sigma)$ be an annihilating pair.
Let $v_0\in L^2(\R)$ and assume that $\widehat{v_0}$ is supported in $\Sigma$.
Let $v$ be a solution of (\ref{eq:shr}) with initial data $v_0$ and assume that,
for some $t>0$, $v(\cdot,t)$ is supported in $S$, then $v_0=0$.

Moreover, if $(S,\Sigma)$ is a strong annihilating pair with constant $D(S,\Sigma)$
in (\ref{eq:sa2}), then, for all $t>0$,
$$
\norm{v(x,t)}_{L^2(\R\setminus S)}\geq D(S,\Sigma)\norm{v(x,t)}_{L^2(\R)}=D(S,\Sigma)\norm{v_0}_{L^2(\R)}.
$$
For instance, if $\Sigma=[-1,1]$ and $E$ is $\gamma$-thick at scale $a>1$, then
$$
\norm{v(x,t)}_{L^2(E)}\geq\left(\frac{\gamma}{C}\right)^{Ca}\norm{v_0}_{L^2(\R)}.
$$


\begin{thebibliography}{99}
\bibitem[AB]{AB}
\textsc{W. O. Amrein \& A. M. Berthier}
\newblock{\em On support properties of $L\sp{p}$-functions and their Fourier transforms.}
J. Functional Analysis {\bf  24} (1977),  258--267.

\bibitem[Ba]{Ba}
\textsc{R. Balian}
\newblock{\em Un principe d'incertitude fort en th\'eorie du signal ou en m\'ecanique quantique.}
C. R. Acad. Sci. Paris S\'er. II M\'ec. Phys. Chim. Sci. Univers Sci. Terre {\bf 292} (1981), 1357--1362.

\bibitem[Be]{Be}
\textsc{M. Benedicks}
\newblock{\em On Fourier transforms of functions supported on sets of finite Lebesgue measure.}
 J. Math. Anal. Appl. {\bf 106} (1985), 180--183.

\bibitem[BD]{BD}
\textsc{A. Bonami \& B. Demange}
\newblock{\em A survey on the uncertainty principle for quadratic forms}
\`a parrai\^itre dans Collecteana Math..

\bibitem[BDJ]{BDJ}
\textsc{A. Bonami, B. Demange \& Ph. Jaming}
\newblock{\em Hermite functions and uncertainty principles for the Fourier and the windowed Fourier transforms.}
Rev. Mat. Iberoamericana {\bf 19}  (2003),  23--55

\bibitem[De]{De}
\textsc{B. Demange}
\newblock{\em Principes d'incertitude associ\'es \`a des formes quadratiques non d\'eg\'en\'er\'ees}
Th\`ese de l'universit\'e d'Orl\'eans, 2004.

\bibitem[DGS]{DGS}
\textsc{P. Delsarte, J.~M. Goethals \& J.~J. Seidel}
\newblock{\em Spherical codes and designs}
Geometrica Dedicata {\bf 6} (1977) 363--388.

\bibitem[EKV]{EKV}
\textsc{L. Escauriaza, C. E. Kenig \& L. Vega}
\newblock{\em Decay at infinity of caloric functions within characteristic hyperplane.}
preprint.

\bibitem[EKPV]{EKPV}
\textsc{L. Escauriaza, C. E. Kenig, G. Ponce \& L. Vega}
\newblock{\em On unique continuation for Schroedinger equation.}
to appear in Comm. Part. Diff. Eq.

\bibitem[Fa]{Fa}
\textsc{W.G. Faris}
\newblock{\em Inequalities and uncertainty principles.}
J. Math. Phys. {\bf 19} (1978), 461--466.

\bibitem[FS]{FS}
\textsc{G. B. Folland \& A. Sitaram}
\newblock{\em The uncertainty principle --- a mathematical survey.}
J. Fourier Anal. Appl. {\bf 3} (1997), 207--238.

\bibitem[Ha]{Ha}
\textsc{G. H. Hardy}
\newblock{\em A theorem concerning Fourier transforms.}
J. London Math. Soc. {\bf 8} (1933), 227--231.

\bibitem[HJ]{HJ}
\textsc{V. Havin \& B. J\"oricke}
\newblock{\em The uncertainty principle in harmonic analysis.}
Springer-Verlag, Berlin, 1994.

\bibitem[JP]{JP}
\textsc{Ph. Jaming \& A. Powell}
\newblock{\em Uncertainty principles for orthonormal bases.}
preprint 2006.

\bibitem[Ko]{Ko}
\textsc{O. Kovrijkine}
\newblock{\em Some results related to the Logvinenko-Sereda theorem.}
Proc. Amer. Math. Soc. {\bf 129} (2001) 3037--3047.

\bibitem[Lo]{Lo}
\textsc{F. Low}
\newblock{\em Complete sets of wave packets.}
in  A Passion for Physics---Essays in Honor of Geoffrey Chew, edited by C. DeTar  et al. (World Scientific, Singapore, 1985), pp. 17--22.

\bibitem[LP1]{LP1}
\textsc{H. J. Landau \& H. O. Pollak,}
\newblock{\em Prolate spheroidal wave functions, Fourier analysis and uncertainty. II.}
Bell System Tech. J. {\bf  40} (1961), 65--84.

\bibitem[LP2]{LP2}
\textsc{H. J. Landau \& H. O. Pollak,}
\newblock{\em Prolate spheroidal wave functions, Fourier analysis and uncertainty. III.
The dimension of the space of essentially time- and band-limited signals.}
Bell System Tech. J. {\bf  41} (1962), 1295--1336.

\bibitem[LS]{LS}
\textsc{V. N. Logvinenko \& Ju. F. Sereda,}
\newblock{\em Equivalent norms in spaces of entire functions of exponential type. (Russian)}
Teor. Funkci\u\i Funkcional. Anal. i Prilo\v zen. Vyp. {\bf 20} (1974), 102--111, 175.

\bibitem[Na]{Na}
\textsc{F. L. Nazarov}
\newblock{\em Local estimates for exponential polynomials and their applications to inequalities of the uncertainty principle type. (Russian)}
Algebra i Analiz  {\bf 5}  (1993), 3--66;  translation in  St. Petersburg Math. J.  {\bf 5}  (1994), 663--717.

\bibitem[Pa1]{Pa1}
\textsc{B. Paneah}
\newblock{\em Support-dependent weighted norm estimates for Fourier transforms.}
J. Math. Anal. Appl.  {\bf 189}  (1995), 552--574.

\bibitem[Pa2]{Pa2}
\textsc{B. Paneah}
\newblock{\em Support-dependent weighted norm estimates for Fourier transforms. II.}
 Duke Math. J.  {\bf 92}  (1998), 35--353.

\bibitem[Po]{Po}
\textsc{A. M. Powell}
\newblock{\em Time-frequency mean and variance sequences of orthonormal bases.}
J. Fourier Anal. Appl. {\bf 11} (2005), 375--387.

\bibitem[Pr]{Pr}
\textsc{J. F. Price}
\newblock{\em Inequalities and local uncertainty principles.}
J. Math. Phys. {\bf 24} (1983), 1711-1714.

\bibitem[PS]{PS}
\textsc{J. F. Price \& A. Sitaram}
\newblock{\em Local uncertainty inequalities for locally compact groups.}
Trans. Amer. Math. Soc. {\bf 308} (1988), 105-114.

\bibitem[Sh1]{Sh1}
\textsc{H.~S. Shapiro}
\newblock{\em Uncertainty principles for basis in $L^2(\R)$},
unpublished manuscript (1991).

\bibitem[Sh2]{Sh2}
\textsc{H.~S. Shapiro}
\newblock{\em Uncertainty principles for basis in $L^2(\R)$},
Proceedings of the conference on {\em Harmonic Analysis and Number Theory},
CIRM, Marseille-Luminy, october 16-21, 2005, L. Habsieger, A. Plagne
\& B. Saffari (Eds). {\it In preparation}

\bibitem[SP]{SP}
\textsc{D. Slepian \& H. O. Pollak}
\newblock{\em Prolate spheroidal wave functions, Fourier analysis and uncertainty. I.}
 Bell System Tech. J.  {\bf 40}  (1961), 43--63. 
 
\bibitem[S1]{S1}
\textsc{D. Slepian}
\newblock{\em Prolate spheroidal wave functions, Fourier analysis and uncertainity. IV.
Extensions to many dimensions; generalized prolate spheroidal functions.}
Bell System Tech. J. {\bf 43} (1964), 3009--3057.

\bibitem[S2]{S2}
\textsc{D. Slepian}
\newblock{\em  Some comments on Fourier analysis, uncertainty and modeling. }
SIAM Rev. {\bf 25} (1983), 379--393.


\end{thebibliography}
\end{document}